\documentclass[11pt]{scrartcl}

\usepackage[english]{babel}
\usepackage[utf8]{inputenc}		

\usepackage{url}

\newcommand{\isdef}{\ensuremath{\mathrel{\mathop:}=}}		

\usepackage{enumerate}
\usepackage{amsfonts, amsmath, amssymb}
\usepackage{amsthm}			

\usepackage{pgf, tikz}
\usepackage{tikz-cd}
\usetikzlibrary{calc}
\usetikzlibrary{arrows, arrows.meta}
\usetikzlibrary{positioning}

\usepackage{wrapfig}
\usepackage{nicefrac}					
\usepackage{color}

\DeclareMathOperator{\Vol}{Vol}

\DeclareMathOperator{\inte}{int}

\newcommand{\menge}[1]{\ensuremath{\mathbb{#1}}}
\newcommand{\N}{\menge{N}}
\newcommand{\Z}{\menge{Z}}
\newcommand{\Q}{\menge{Q}}
\newcommand{\R}{\menge{R}}

\newtheorem{lem}{Lemma}[section]		
\newtheorem{ex}[lem]{Example}		

\newtheorem{theo}[lem]{Theorem}			
{\begin{proof}[Well--defined]}
	{\end{proof}}

\newcommand{\ie}{i.\,e.,}

\begin{document}
\noindent\LARGE{Effective Khovanskii, Ehrhart Polytopes, and the Erd\H{o}s Multiplication Table Problem}
\vspace{0.5cm}

\large

Anna M. Limbach$^{1}$,\hspace{.2cm}
Robert Scheidweiler$^2$,\hspace{.2cm}
Eberhard Triesch$^3$
\vspace{0.5cm}
\normalsize

\noindent$^1$ Institute of Computer Science, Czech Academy of Sciences, Pod Vodárenskou věží 271/2, 182 00 Praha 8, Czech Republic\\
$^2$ Faculty of Electrical Engineering \& Information Technology, 
Hochschule Düsseldorf - University of Applied Sciences, Münsterstraße 156, 40476 Düsseldorf, Germany\\
$^3$ Lehrstuhl II für Mathematik, RWTH Aachen University, Pontdriesch 10-12, 52062 Aachen, Germany\\

\section{Abstract}
	Let $P(k,n)$ be the set of products of $k$ factors from the set $\{1,\ldots , n\}.$ In 1955, Erd\H{o}s 
	posed the problem of determining the order of magnitude of $|P (2, n)|$ and
	proved that $|P (2, n)| = o(n^2 )$ for $n \to\infty$.  
	In 2015, Darda and Hujdurović 
	asked whether, for each fixed $n$, $|P (k, n)|$ is a polynomial
	in $k$ of degree $\pi(n)$ - the number of primes not larger than $n$.
	Recently, Granville, Smith and Walker published an effective version of Khovanskii's Theorem.
	We apply this new result to show, that for each integer $n$, there is a polynomial $q_n$ of degree $\pi(n)$ such that $|P (k, n)|=q_n(k)$ for each $k\geq n^2\cdot\left(\prod_{m=1}^{\pi(n)} \log_{p_m}(n)\right)-n+1.$ Moreover, we give an upper estimate of the leading coefficient of $q_n$.

	Keywords: Khovanskii's theorem, Ehrhart polynomial, multiplication table problem, additive combinatorics, integral polytopes
	
\section{Introduction}

	Let 	\begin{align*}
	P(k,n)
	&=\left.\left\{\prod_{j=1}^{n}j^{\lambda_j}\right\vert \lambda_j\in\N_0,\sum_{j=1}^{n}\lambda_j=k \right\}
\end{align*}
	be the set of products of $k$ factors from the set $\{1,\ldots , n\}$ and denote by
$p(k, n)$ its cardinality. In 1955 Erdős considered the problem of estimating
$p(2, n)$ and showed that $p(2, n) = o(n^2 )$ (cf. \cite{erdHos1955some}). His estimates for this so-called Erdős Multiplication Table Problem were considerably refined in 2008 by
Ford in \cite{ford2008distribution} and, for general fixed $k$, by Koukoulopoulos \cite{koukoulopoulos2010localized}. Motivated by a
graph coloring problem, namely determining the product irregularity strength
of complete bipartite graphs, Darda and Hujdurović  asked in \cite{darda2015bounds} whether, for fixed
$n$, the quantity $p(k, n)$ is a polynomial in $k$ of degree $\pi(n)$ - the number of primes not larger
than $n$. They showed this to be true for $n = 1, . . . ,10$. 
We apply Khovanskii's Theorem and Ehrhart's Theorem to answer this question affirmatively for each $n$ and for every $k$ larger than some threshold $k_n$, which depends on $n$. Moreover, we give an expression as well as an upper bound for the polynomials leading coefficient.
By applying a recent result from  Granville, Smith and Walker, we establish an upper bound on $k_n$.

\section{Reformulation by Prime Factorization}

	Let $p=(p_1,\ldots,p_{\pi(n)})$ be the $\pi(n)$-tuple of prime numbers at most $n$.
	We define
	\begin{align*}
		M_n&\isdef \left\{\alpha\in\N_0^{\pi(n)}\left\vert \prod_{j=1}^{\pi(n)}p_j^{\alpha_j}\leq n\right.\right\}\\
	\end{align*}
	to be the set of $\pi(n)$-tuples of $p$-adic evaluations of numbers not exceeding $n$.
	Given that $\{1,\ldots,n\}=\{\prod_{j=1}^{\pi(n)}p_j^{\alpha_j}\mid \alpha\in M_n\}$ and that multiplication can be expressed as addition of exponents, it follows that $$P(k,n)=\left\{\left.\prod_{j=1}^{\pi(n)}p_j^{\alpha_j}\right\vert \alpha\in kM_n\right\},$$ where we use the notation $$kA := \{a_1 + \ldots + a_k : a_i \in A \text{ for all }i\}$$ for a set $A$ and a positive integer $k$.
	This implies $p(k,n)=\lvert P(k,n)\rvert=\lvert kM_n\rvert.$
	
	However, the size of sumsets of the form $kA$ is well understood.
	A famous theorem by Khovanskii  says the following:
	
	\begin{theo}[\cite{khovanskii1992newton}]
		Let $A \subseteq \Z^d$ be finite. There is a polynomial $P_A \in \Q[X]$ of
		degree at most $d$, and a threshold $N_{Kh}(A)$, such that $|NA| = P_A(N)$ provided $N > N_{Kh}(A)$.
	\end{theo}

	Thus we obtain the following intermediate result:
	
	\begin{lem}
		For each $n\in \N$ there is a polynomial $q_n\in\Q[X]$ of degree at most $\pi(n)$ and a threshold $k_n$ such that $p(k,n)=q_n(k)$ for each $k\geq k_n.$
	\end{lem}

\section{Relation to Ehrhart Polynomials}

We now consider a \emph{full dimensional lattice polytope} $Q$ in $\R^d$, \ie\ a polytope whose vertices are
all integral and which contains $d$ linearly independent vectors. Denote by $\inte(Q)$ the set of integral points in $Q$. Expanding
$Q$ by a factor of $t \in \N$ in each dimension, the number of integral points in $tQ$
is denoted by $L(Q, t) = |\inte(tQ)|$. As was proven by Ehrhart in 1962 \cite{ehrhart1962polyedres}, the
function $L(Q, t)$ is a polynomial in $t$ of degree $d$, the Ehrhart polynomial of $Q$.
It is well-known that the coefficient of $t^d$ in $L(Q, t)$ is the $d$-dimensional volume
of $Q$. Denote by $k \star Q=k(\inte(Q_n))$ the set of all vectors in $\R^d$ which can be obtained as
a sum of $k$ integral vectors in $Q$. Clearly, $k \star Q \subseteq kQ$ and the polytope $Q$
is called \emph{integrally closed} if $k \star Q = \inte(kQ)$ for all $k$. There are very simple
polytopes which are not integrally closed. The following example is well-known
(\cite{bruns2009polytopes}, Example 2.56 (c)):
\begin{ex}
	Let $Q := \{(x, y, z) \in\R^3 : x \geq 0, y \geq 0, z \geq 0, 6x + 10y + 15z \leq 30\}$.
	As is easily verified, the integral points in $Q$ are $$(5, 0, 0), (3, 1, 0), (2, 0, 1),
	(1, 2, 0), (0, 3, 0), (0, 1, 1), (0, 0, 2)$$ and component-wise smaller points in $\N^3_0$. The
	point $(4, 2, 1)=2\cdot(2,1,1/2)$ is in $2Q$, but not in $2 \star Q$, hence the simplex $Q$ is not integrally
	closed.
\end{ex} 
Now let $Q_n$ denote the convex hull $H(M_n)$, which is a polytope of dimension $\pi(n)$. Then $k \star
Q_n = kM_n$ and thus $p(k,n) = |kM_n| \leq L(Q_n,k)$. In fact, we checked
that $Q_n$ is integrally closed and hence the polynomials $L(Q_n,t)$ and $q_n(t)$ are
equal for $n \leq 20$. By adapting the example above, we can, however, show that
the polytopes $Q_n$ are not integrally closed in general:
\begin{ex}
 Choose some large prime $p$ and let $n = p^5$. By Bertrand’s postulate
(see \cite{chebyshev1852memoire}, easier proofs can be found in \cite{ramanujan1919proof} and \cite{erdos1932beweis}), we can choose primes $q$ and $r$ such
that $n^{1/3} \leq q \leq 2n^{1/3}$ and $n^{1/2} \leq r \leq 2n^{1/2}$. We write the points in $Q_n$ in the
form $(x, y, z, . . .)$ where the coordinates $x, y, z$ correspond to the primes $p, q, r,$
respectively. Then a nonnegative integral vector $(x, y, z, 0, . . . , 0)$ is in $Q_n$ if and
only if
$$\frac{\log(p)}{\log(n)}x+\frac{\log(q)}{\log(n)}y+\frac{\log(r)}{\log(n)}z\leq 1.$$
This is equivalent to

$$\frac15x+(\frac13+\epsilon)x+(\frac12+\delta)z\leq 1$$
where $0 \leq \epsilon,\delta \leq \log(2)/\log(n) < 1/60$ for large $n$. As in the example above,
$(4, 2, 1, 0, . . .,0)$ is in $2Q_n$, but not in $2M_n$.
\end{ex}

We now show that the coefficients of $t^d$ in $q_n(t)$ and $L(Q_n,t)$ are equal:

\begin{theo}
	\begin{enumerate}[(i)]
		\item  For all $k,n\in \N$, the inequalities
		$$p(k, n) \leq L(Q_n,k) \leq p(k + \pi(n), n)$$
		hold.
		\item The polynomial $q_n$ has degree $\pi(n)$ with leading coefficient equal to
		the $\pi(n)$-dimensional volume of the polytope $Q_n$.
	\end{enumerate}
\end{theo}

\begin{proof}
		\begin{enumerate}[(i)]
		\item  Let $d\isdef \pi(n)$. We show that $kM_n \subseteq \inte(kQ_n) \subseteq (k+d)M_n$, where the
		first inclusion was already explained above. For the second, recall that, by
		Caratheodory’s Theorem, each point $x \in \inte(kQ_n)$ can be written as a linear combination $x = \sum_{i=1}^d \alpha_i v_i$ with $v_i \in M_n$, $\alpha_i \geq 0$
		for all $i$ and $\sum_{i=1}^d \alpha_i = k$. It follows that $\tilde{x} := \sum_{i=1}^d \lceil\alpha_i\rceil v_i$ is in $(k+d)M_n$. The set $M_n$ is a down-set, \ie\ for each $y\in M_n$ and each $z\leq y$ (component-wise), we obtain $z\in M_n$. Consequently, $kM_n$ and $(k+d)M_n$ are down-sets as well and $x \leq \tilde{x}$ implies that $x$ is in $(k+d)M_n$ as well.
		\item This follows immediately from (i) and from the
properties of Ehrhart polynomials.
	\end{enumerate}

\end{proof}

\section{Application of the Effective Khovanskii Theorem}

In a new preprint, Granville, Smith and Walker prove an effective version of Khovanskii's theorem, \ie\ an upper bound on $N_{Kh}$:

\begin{theo}[\cite{granville2024improved}]
	Let $A \subseteq \Z^d$ be finite and let $H(A)\subseteq \R^d$ be its convex	hull. Then, $$N_{Kh}(A)\leq d!|A|^2 \Vol(H(A))-|A| + 1.$$
\end{theo}
By the injectivity of the logarithm on $\R_{\geq 0}$ we obtain $$M_n=\left\{\alpha\in\N_0^{\pi(n)}\left\vert\ \sum_{j=1}^{\pi(n)}\alpha_j\cdot \log(p_j)\leq \log(n)\right.\right\},$$
thus, it is easy to see that $$H(M_n)\subseteq \left\{x\in \R_{\geq 0}^{\pi(n)}\left\vert\ \sum_{j=1}^{\pi(n)}x_j\cdot \log(p_j)\leq \log(n) \right.\right\}$$ and hence $\Vol(H(M_n))\leq \frac{1}{\pi(n)!}\prod_{j=1}^{\pi(n)}\log_{p_j}(n)$.\\[0.2cm]
We obtain the following main result:

\begin{theo}
	For each $n\in \N$ there is a polynomial $q_n\in\Q[t]$ of degree $\pi(n)$ with leading coefficient $\Vol(H(M_n))$ such that $p(k,n)=q_n(k)$ for each $k\geq n^2\cdot\left(\prod_{m=1}^{\pi(n)} \log_{p_m}(n)\right)-n+1.$ Furthermore, the leading coefficient fulfills $ \Vol(H(M_n))\leq \frac{1}{\pi(n)!}\prod_{j=1}^{\pi(n)}\log_{p_j}(n)$. 
\end{theo}

\section{Future Research}

The result by Granville, Smith, and Walker is quite general, as it can be applied to every finite subset of $\Z^d$.
We wonder, whether one can get tighter bounds by modifying their proof and using more properties of $M_n$, for example that $M_n$ is downward-closed or that $M_n$ is the set of integer vectors of a polytope.

%
%

\bibliography{referenzen}
\bibliographystyle{plain}

\end{document}